\documentclass[11pt]{article}
\usepackage{amsmath,amssymb}
\usepackage{enumerate}
\usepackage{color}

\usepackage{xy}
\usepackage{xypic}
\usepackage{subfigure}
\usepackage{tikz-cd}

\newcommand{\R}{\mathbb{R}}

\newcommand{\N}{\mathbb{N}}
\newcommand{\bea}{\begin{eqnarray}}
\newcommand{\eea}{\end{eqnarray}}
\newcommand{\la}{\lambda}

\def\a{\alpha}

\def\de{\delta}
\def\e{\varepsilon}

\def\s{\sigma}

\def\diam{{\rm diam}}

\def\1{\rm Id}

\def\sup{{\rm sup}}
\def\vol{{\rm vol}}
\def\Vol{{\rm Vol}}

\def\ci{\circ}
\def\cl{{\rm cl}}

\newcommand{\qed}{$\hfill\blacksquare$}

\def\V{\noindent}

\def\csec{{\rm csec}}
\def\cdiam{{\rm cdiam}}

\def\injrad{{\rm injrad}}
\def\Lor{{\rm Lor}}

\def\SEP{{\rm sep}}
\def\F{\rm{fld}}
\def\KS{{\rm fks}}
\def\MCS{\rm{mcs}}

\def\cl{{\rm cl}}
\def\Dm{{\rm Dm}}

\def\CPOM{\rm {\bf POM}}

\def\CGH{\rm {\bf GHM}}
\def\CMet{\rm {\bf MET}}
\def\CTop{\rm {\bf TOP}}
\def\CALP{\rm{\bf ALP}}
\def\CALL{\rm {\bf ALL}}
\def\CCS{\rm {\bf CST}}

\def\CR{{\bf R}}
\def\CC{{\bf C}}

\def\CW{{\bf WLS}}
\def\CWP{{\bf WPL}}

\newcommand{\bean}{\begin{eqnarray*}}
\newcommand{\eean}{\end{eqnarray*}}

\newtheorem{Theorem}{Theorem}

\newtheorem{Definition}{Definition}

\newcommand{\ben}{\begin{enumerate}}
\newcommand{\een}{\end{enumerate}}
\newcommand{\bit}{\begin{itemize}}
\newcommand{\eit}{\end{itemize}}
\newcommand{\edoc}{\end{document}}

\usepackage{hyperref}

\usepackage{enumerate}

\title{Functors in Lorentzian geometry: three variations on a theme}

\begin{document}

\author{Olaf M\"uller\footnote{Institut f\"ur Mathematik, Humboldt-Universit\"at zu Berlin, Unter den Linden 6, D-10099 Berlin, \texttt{Email: mullerol@math.hu-berlin.de}}}

\date{\today}
\maketitle

\begin{center}
	{\em Dedicated to Roger Penrose}
\end{center}

\begin{abstract}
\V We consider three examples of functors from Lorentzian categories and their applications in finiteness results, singularity theorems and boundary constructions. The third example is a novel functor from the category of ordered measure spaces to the category of Lorentzian pre-length spaces in the sense of Kunzinger-S\"amann.    
\end{abstract}

Categories and functors are important structuring elements in geometry and other branches of mathematics and often allow for a concise formulation of complex situations. The laureate of this conference volume has some interesting relations to category theory, two sample examples given by the Penrose transform in twistor theory (\cite{rP:PT1}, \cite{rP:PT2}, \cite{rP:PT3}) and monoidal categories applied to spin networks (see e.g. \cite{rP:SN}). Here, we will not deepen further those topics but rather give some recent advances in the theory of functors with applications related to Roger Penrose's work.

We want to present three functors from Lorentzian categories, the first one mapping into a Riemannian category with applications in finiteness theorems, the second one from ordered sets (encoding a causal structure) to topological sets, with application in the theory of black holes, and the third one from the category of ordered measure spaces to the category of (almost) Lorentzian length spaces, with an application yet to be found. All these functors factorize over the identity in the category of sets. The author would like to apologize in advance for not including many other interesting functorial approaches to Lorentzian geometry, like the one taken in \cite{mG} where the conformal structure of a spacetime is fully encoded by the isomorphism class of the category of causal path homotopies, and conformal maps are functors in this category.

 \V Whereas the first two sections of the article review results stated and proven elsewhere (except for a short paragraph following Th. \ref{Main1}), the theorems in the third chapter are entirely original material.

The author gratefully acknowledges a helpful conversation with Clemens S\"amann and constructive remarks of two anonymous referees on previous versions of the article.

\section{Lorentzian-to-Riemannian functors and finiteness results}
\label{LorFin}

Lorentzian and Riemannian geometry look like two worlds apart. The Riemannian side is dominated by elliptic equations with their smoothing properties, whereas Lorentzian geometry is the natural arena for hyperbolic equations, preserving all kinds of singularities. On the Lorentzian side, we have non-convexity of the space $\Lor (M)$  of Lorentzian metrics on a fixed manifold $M$ within the space of bilinear forms, the notorious difficulty to find appropriate topologies on $\Lor(M)$, the lack of a Hopf-Rinow statement, and the fact that isometric actions might have noncompact isotropy groups and indecomposable yet reducible subspaces. On the other hand, we have a nice link between causality and Lorentzian conformal structure, null pregeodesics are conformally invariant, the initial  value problem for Lorentzian minimal surfaces on globally hyperbolic spacetimes is better treatable than the corresponding Riemannian one (\cite{oM-stringy}) and last not least we have singularity theorems.

M. Berger writes in his book "Panoramic View of Riemannian Geometry" (\cite{mB}, Sec. 14.6):

"{\em Gromov's mm} [= metric measure] {\em spaces are, in our opinion, the geometry of the future.}"

Indeed, widening the scope and considering Riemannian manifolds as special metric spaces has turned out to be a fruitful method in Riemannian geometry e.g. for finiteness theorems, where one can use Gromov compactness \cite{mG} and Perelman stability \cite{Kapovitch}. The natural question arises whether some analogous procedure is possible for Lorentzian geometry as well. To the author's knowledge, the first full-fledged synthetic approach to Lorentzian geometry is the one of Lorentzian length spaces by Kunzinger and S\"amann \cite{mKcS}, which is fascinating and appropriate for many aspects, but likely not of great use for finiteness theorems, as its metrical component is not functorial, and its Lorentzian distance component, which {\em is} functorial, satisfies an {\em inverse} triangle inequality. 

\smallskip

The first discouraging fact in this context is the following folk wisdom:

\begin{equation}
	\label{NOGO}
{\rm There \ is \  no \ functor\ } {\rm F} {\rm \ with} \	{\rm top} \ci {\rm F} = {\rm forget}.
\end{equation}

where ${\bf forget}$ is the forgetful functor from the category $\CGH$ of globally hyperbolic Lorentzian manifolds and time-oriented isometries to the underlying topological spaces, and ${\rm top} $ is the functor from the category $\CMet$ of metric spaces to the category $\CTop$ of topological spaces assigning to each metric the topology it generates. The simple reason is that the orbit of any nonzero null vector by $SO(1,n)$ and thus also by ${\bf F}(SO(1,n))  \subset {\rm Isom} ({\bf F} (\R^{1,n})) $ would accumulate at $0$ (a boost, i.e., an element of $SO(1,1) \subset SO(1,n)$, maps $u \mapsto \lambda u , v \mapsto \lambda^{-1} v$ in null coordinates $u,v$). As each temporal function on a spacetime gives rise to a Riemannian metric, this means that there is no canonical choice of temporal function (but there are canonical {\em classes} of temporal functions with preferable properties, see e.g. \cite{oM-special}). However, we can exclude boosts by replacing $\CGH$ with the category of {\bf Cauchy slabs}, which are globally hyperbolic manifolds whose boundary consists of two disjoint connected smooth spacelike Cauchy surfaces: We define

\bean
\CC_n^-&:=& ( \{  {\rm smooth \ } n {\rm -dimensional \ Cauchy \ slabs} \}, \{{\rm isometric \ oriented \ time-oriented \ diffeomorphisms}\} ),\\
\CC_n^+ &:=& ( \{ {\rm smooth \ } n-{\rm dimensional \ Riemannian \  manifolds-with-boundary } \}  , \{ {\rm oriented \ isometries }\}   ).
\eean

The following results of this section are shown in \cite{oM:LGH}. Let $\s^g: X \times X \rightarrow \R$  the signed Lorentzian distance $\sigma^g : X \times X \rightarrow \R \cup \{ \infty \} \cup \{ - \infty \}$ defined by $\s^g(x,y) := \mp \sup \{ \ell(c)  \vert c: x \leadsto y {\rm \  causal \  curve} \}$ for $x \in J^\pm ( y)$ and $\s^g (x,y) = 0$ otherwise, where $\ell$ is Lorentzian length and $J^\pm (y)$ the causal past resp. future of $y$. We put $\s_x:= \s^g(x, \cdot)$, and we define, for 

$f \in {\rm Adm} := \{ f: \R\rightarrow \R \  {\rm  \ measurable \ locally \ essentially \  bounded, \ } f^{-1} (0) = \{ 0\}, f \vert_{\pm (0; \infty)} {\rm injective} \}:$ 

\bea
\label{KuratowskiEmb}
\Phi_{f}^g: x \mapsto f \ci \sigma_x  \in L^p (X) \ \ \forall p \in [1; \infty] .
\eea

We pull back the extrinsic/intrinsic metrics of Banach spaces on the right-hand side of Eq. \ref{KuratowskiEmb}. For $p :=1$, $f:= \chi_- := \chi_{(- \infty; 0 )} $ (with $\chi_A$ the characteristic function of $A \subset \R$), we get the {\em Beem metric}

$$d_{\chi_- , 1}(p,q):= \vert \chi_- \circ  \sigma_p  - \chi_- \circ \sigma_q \vert_{L^1} = \vol ((I^- (p) \setminus I^-(q) ) \cup (I^-(q) \setminus I^- (p))) ,$$ 

which will also be of central importance in Section 2. It has splitting geodesics.

For $p:=\infty$, $ f:= \vert . \vert $, we get the {\em Noldus metric} \cite{Noldus}, not inducing the manifold topology \cite{oM:LGH}.

Here we choose $p:=2, f \in \{  f_r:= (\frac{1}{2} + \frac{r}{2} {\rm sgn}) \cdot \vert . \vert \cdot \1 ^3 , h := \1^4 \}$. The reason for choosing quartic polynomials is that many derivatives at $0$ vanish, ensuring that $\Phi_g^f$ is $C^2$. It would be desirable to work with $h$ alone, but in order to reconstruct the Lorentzian isometry class from the data on the right-hand side it turns out to be useful to include the $f_r$ as additional data and define an enriched category $\CR_n$ of objects $X$ of $C_n^+$ with Lipschitz functions $X \times X \rightarrow \R^3 $, whose isometric morphisms are additionally required to pull back the functions. We obtain, with $d_r:= (\Phi_{f_r}^g)^* d_{L^2(X)}$:

\begin{Theorem}[\cite{oM:LGH}]
	Let $X$ be a Cauchy slab.
	
	\begin{itemize}
		
		\item $ \forall r \in (-1; 1): \Phi_{f_r}^g$ and $\Phi_h^g$ are $C^2$ embeddings $X \rightarrow L^2(X)$.
		
		\item $ F: (X,g) \mapsto (X, (\Phi_{h}^g)^* (\langle \cdot, \cdot \rangle_{L^2 (X)}), (d_{-1/2}, d_0, d_{1/2}) )$ is an injective functor $ {\mathbf C}_n^- \rightarrow \CR_n$.

	\end{itemize}	
	
\end{Theorem}

Let $X \in \bf{Obj} (C_n)^\pm$, $\nu $ be the outer normal at $\partial X$, $\csec$ := infimum of cospacelike sectional curvature (for $X \in \CC_n^+$ all planes are cospacelike, s.t. $\csec$ is then just the infimum of sectional curvature),  

$\cdiam (X) :=  
\begin{cases}
	 \diam (X) ,& X \in \CC_n^+ \\
	\sup \{ g(\nu (p) , w (p)) \vert p \in \partial X, w \in I_p, \exp_p (w) {\rm \ exists } \} \geq \sup \{  \s(x,y) \vert x,y \in X \}   ,& X \in \CC_n^- ,
\end{cases}$

\V where $I_p$ is the set of timelike vectors in $T_pM$. (For the last inequality let $D := \cdiam (X)  $, we show $\de:=  \sup \{ \s (x,y) \vert x,y \in X \} < D$: the inverse triangle inequality implies that we can restrict the supremum to distances from points in $\partial^- X$ to points in $\partial^+ X$. For given $\e >0$ we find $x_\pm \in \partial^\pm X$ with $d(p_-, p_+) > \de - \e$. By compactness of $J^-(p_+)$ we find $\tilde{p}_- \in \partial^- X$ with $d(\tilde{p}_-, p_+ ) = \sup \{ d(x,p_+) \vert x \in X \} > d(p_-, p_+) > \de - \e$. Let $c: [0;1] \rightarrow X$ be a maximizing curve from $\tilde{p}_-$ to $p_+$. By the first variational formula, there is $k \in (0; \infty) $ with $c'(0) = k \cdot \nu (c(0))$, thus 

$\forall \e > 0 : \de - \e < d(\tilde{p}_-, q) = \ell (c) = \sqrt{- g(c'(0), c'(0))}= k = g (c'(0), \nu (c(0))) \leq D $, so $\de \leq D$.)

We define ${\bf C}_n^{+,{\rm c}}  $ to be the subcategory of metrically complete objects in ${\bf C}_n^+$, and ${\bf C}_n^{-,{\rm c}} := {\bf C}_n^-$,

$J(p) := J^+(p) \cup J^-(p)$, \ ${\rm JV} (X) := \sup \{ \vol (J (p) ) \vert p \in X \}, \ \ {\rm JV} (\partial X) := \sup \{ \vol (J (p) \cap \partial X  ) \vert p \in X \}  $,

\bigskip

	${\rm injrad}^\pm_g (x) :=  \sup \{  \sqrt{- g(w,w)} \big\vert \  w \in T_xX, w \ll 0, \exp\vert_ {I^\pm(0) \cap I^\mp(w)} {\rm \ is \ a \   diffeomorphism } \}  , $
	
\medskip

$ {\rm injrad} (X) := \inf_{x \in X} \max \{ {\rm injrad}^+_g (x) , {\rm injrad}^-_g (x)   \} .$ 

\medskip

Then $g \mapsto \Phi^g$ transfers upper bounds on the above data: For $a \in \R^4$, $b \in  \R^7$ we define

\medskip
	
$	C_n^\pm (a) := \{ (X,g) \in {\rm Obj} (C_n^{\pm, c}) \big\vert  - \csec (X) \leq a_1, \vert \nabla \nu \vert_X \leq e^{a_2},
	\cdiam (X) \leq e^{a_3},  ({\vol} (X))^{-1} \leq e^{a_4} \} ,$
	
$ C_n^- (b) := \{ (X,g) \in C_n^- (b_1,b_2, b_3, b_4)   \big\vert  {\rm JV}^g (X) \leq e^{b_5}, \injrad (X) \geq e^{-b_6} ,
	{\rm JV}^g ({\partial} X ) \leq e^{b_7}  \}.$	

\medskip

The considered conditions do {\em not} imply, by mere restriction, estimates on the Cauchy surfaces that would allow for an application of Riemannian finiteness theorems for them. However, we obtain:

\begin{Theorem}[\cite{oM:LGH}]
	\label{Main1}
	Let $n \in \N $ and $X \in {\rm Obj} ({\bf C_n^-})$.
	\begin{enumerate}
		\item (Bound transfer) $  \forall \  b \in \R^7: \big( b_1 \geq 0 \Rightarrow \exists \ a \in \R^4 \   \forall (X,g) \in C_n^-(b) : \ 
		(X, (\Phi^g)^* \langle \cdot , \cdot \rangle ) \in C_n^+ (a) \big) $. 
		\item (Richness) $ \forall \  b \in \R^7  : b_6 \leq b_3 \land b_4 \leq b_5 \Rightarrow C_n^-(b) $ contains an open set in $C_n^-$. 
	\end{enumerate}
\end{Theorem} 

In other words, $\Phi$ transfers an appropriate bound $b$ from the Lorentzian category to a bound $a$ on the Riemannian category. Let us be more explicit than \cite{oM:LGH} about the completeness of the induced Riemannian metrics: It suffices to show that any Cauchy sequence $a$ in $X$ remains in a compact subset of $X$. This latter assertion follows from the fact that eventually any two $J(a(n))$ have to intersect each other, which can be seen by estimating $\int \vert \s(a(n), x) \vert dx $ from $0$ by the bounds on $\injrad^-_g$ and $\csec$ exactly as in the proof of Lemma 3 of \cite{oM:LGH}.

\begin{Theorem}[\cite{oM:LGH}]
	\label{Main2}
	For every $n \in \N$ and for all $ b \in \R^7$ with $b_1 \geq 0$, there are only finitely many homeomorphism classes of compact Cauchy slabs in $C_n^-(b)$.
\end{Theorem}

The article also obtains auxiliary results about the dependence between the second fundamental form, the intrinsic and the extrinsic diameter of a submanifold-with-boundary $M$ in a Hilbert space $H$ that are of independent interest. For the special case $\partial M = \emptyset $ we get e.g.

\begin{Theorem}[\cite{oM:LGH}]
	\label{Hilbert}
	Let $H$ be a Hilbert space. Let $M$ be a closed (boundaryless) submanifold of $H$ with second fundamental form $S_M^H$ and extrinsic resp. intrinsic diameter $\de$ resp. $D$. \\ If $\vert \vert S_M^H \vert \vert < \frac{4 \sqrt{2}}{3 \de}$, then $D < 3 \de$.  
\end{Theorem}

The dependencies of upper bounds are listed in the table below (e.g. $\csec (X,g) \geq 0$ and upper bounds on on $ 1/\vol $ and $1/\injrad$ of $(X,g)$ imply an upper bound on $1/\vol$ of $(X, (\Phi^g)^* \langle \cdot , \cdot \rangle )$): 

\bigskip


\renewcommand{\arraystretch}{1.1}
\begin{tabular}[h]{c||c|c|c|c|c|c|c||}
	\  & -csec $\leq$ 0 & $\vert \nabla \nu \vert $ & 1/vol & cdiam & Jvol & 1/injrad & Jvol $\circ \partial$\\
	\hline 
	\hline
	-csec $\leq$ 0 & x & \  & \  & \ & \ & \ & \  \\
	\hline
	$\vert\nabla  \nu \vert  $ & \ & x & x & \ & x & \ & x \\
	\hline
	1/vol & x & \  & x & \ & \ & x & \ \\
	\hline
	cdiam & x & x & x & x & \ & x & \  \\
	\hline
	\hline
\end{tabular}
\renewcommand{\arraystretch}{1}

\newpage

\section{Functors from orders to topologies, applied to black hole theory}

Usually, future null infinity and black holes are defined via asymptotic simplicity, which assumes the existence of a conformal boundary of spacetimes. Unfortunately, the latter is an often unfounded assumption --- but let us see how far we get heading for nice conformal boundaries. Maximally, we would ask for a generalization of the Penrose embedding (which is the canonical extension of stereographic projection of a linear slice to the entire Minkowski spacetime): an open conformal embedding $F$ of the spacetime $X$ into a larger globally hyperbolic spacetime whose image is relatively compact (and, if possible, even causally convex). The rigidity statement of the positive mass theorem tells us that such an embedding cannot exist if $X$ is any maximal Einstein-Maxwell Cauchy development of an initial value with nonvanishing electromagnetic field --- positivity of mass is an obstruction to smoothness of the conformal structure at spatial infinity (see e.g. \cite{oM:BH}). The first concession we could make is that we replace compactness of the image of $F$ by the more flexible one of strong future compactness, which requires only that the image of $F$ is contained in the past of a compact set. Such a map exists for every maximal Cauchy development of initial values of Einstein-Maxwell theory small in a certain weighted Sobolev norm, and still, this condition implies global existence of massless Dirac-Higgs-Yang-Mills systems for every initial value small in weighted Sobolev norms \cite{nGoM} and leads to existence results on black holes in Einstein-Maxwell theory \cite{oM:BH}. Nevertheless, even this more liberal definition is not always satisfied --- one example showing this defect is Kruskal spacetime, the simplest example of a black hole. Another, actually quite old, idea to overcome this deficiency is to define an {\em abstract}  future boundary of spacetimes, called the Geroch-Kronheimer-Penrose boundary \cite{GKP} $IP(X) \setminus X$. $IP(X)$ consists of all {\em indecomposable past sets}, that is, all past sets that cannot be written nontrivially as a union of two past sets, and the boundary consists of all such subsets that are not the past of a point of $X$. Geroch, Kronheimer and Penrose \cite{GKP} found that they are always pasts of continuously inextendible future timelike curves. There is an obvious canonical causal structure on $IP(X)$ defined by $A \leq B : \Leftrightarrow A \subset B$, and a canonical chronological relation $\ll$ due to Budic and Sachs \cite{BS} defined by $  A \ll B : \Leftrightarrow  \exists p \in B: A \subset I^-(p) $. 

We can also induce a causal relation $\a  (\ll)$ from a chronological relation $\ll$, by $x \a (\ll) y: \Leftrightarrow I^-(x) \subset I^-(y) $. In \cite{oM:TFC}, we show that on $IP(X)$ both notions of inducing a causal structure coincide for $X$ causally simple (and otherwise, there can be spurios causal relations, as in the example $\R^{1,n} \setminus \{ 0 \} $).

The converse perspective is to induce a chronological relation $\beta (\leq)$ from a causal relation (a partial order) $\leq$, e.g. following Minguzzi and S\'anchez (\cite{eMmS}) via $ p <q : \Leftrightarrow p \neq q \land p \leq q $ and

$(x,y) \in \beta (\leq) : \Leftrightarrow \big( x \leq y \land ( \exists u,v \in X:  x < u< v < y \land J^+(u) \cap J^-(v) {\rm \ not \ totally \ ordered} ) \big).$

We could topologize $IP(X)$ by the Alexandrov (interval) topology, having chronological diamonds $I^+(x) \cap I^-(y)$ as a subbasis. Unfortunately, this would deviate from the topology obtained by conformal completions, even if we include as neighborhoods $I^+(z)$ and/or causal diamonds $J^+(x) \cap J^-(y) $, as it can be seen in the example of $\R^{1,1} \setminus I^+(0) = IP(\R^{1,1} \setminus J^+(0))$,  where $0$ would not have enough neighborhoods. A better proposal is the {\em chronological topology} (due to Flores inspired by Harris, \cite{jF}). 
It is defined via
$C \subset IP(X) $ being $\tau_-$-closed if and only if

$$ \forall a: \N\rightarrow C:  L_- (a)  := \{ P \in IP(M) \setminus \{ \emptyset \}  \big\vert P \subset \liminf (a) \land P {\rm \ maximal \   in \ } IP(\limsup (a) )\} \subset C.$$

(with set-theoretic liminf and limsup). One of its advantages is that it is very generally applicable: For all strongly causal spacetimes $X$, the canonical injection $i_X : X \rightarrow (IP(X), \tau_-)$ is a homeomorphism on its dense image. It is very well compatible with the full abstract completion (which includes the past boundary with some canonical identification between future and past boundary points). The main drawback is that $(IP(X), \tau_- ) $ is non-Hausdorff in general, even for g.h. spacetimes, e.g. for $-dt^2 + g$ where $g$ is the metric on the unwrapped grapefruit-on-a-stick (an example due to S. Harris, \cite{sH}). 

Another proposal, which is actually older but apparently had been forgotten by the community, has been made by Beem \cite{Beem}. It uses a distance $d_1$ on the space $C(X)$ of closed subsets of a metric space $(X,\de) $ that had been defined by Busemann via

\bea
\label{LengthDefTau1}
d_1(A,B) := \sup \{  \vert \de ( \{ x \},A) - \de (\{ x \},B) \vert \cdot \exp (-\de( x_0 , x)) : x \in X \}  \ \forall A,B \in C(X),
\eea

for some $x_0$ on which the equivalence class of $d_1$ does not depend. If $\de$ is Heine-Borel, so is $d_1$. We define a metric on $IP(X)$ by $d(A,B):= d_1(\cl A , \cl B)$ and the induced topology $\tau_+$, which can also be defined via a finite regular Borel measure $\mu$ on $X$ and the symmetric difference $\triangle$ of sets by the Beem metric (w.r.t. the measure $\mu$), which assigns to each two $p,q \in X$ the number

\bea
\label{MeasureDefTau1}
 d(p,q) := d_{\chi_{(- \infty, 0)},1}(p,q)  = \Vol_\mu (I^-(p) \triangle I^-(q) ) .
 \eea

In summary, a topology on $IP(X)$ can be induced by using either a distance on $X$ or a finite measure on $X$. Now, there is no natural choice of neither distance (see Eq. \ref{NOGO}) nor measure of finite mass on $X$ (however, compare the last section), but one can show that the topology $\tau_+$ induced by Eq. \ref{LengthDefTau1} or \ref{MeasureDefTau1} is independent of the choices and that $C \subset IP(X) $ is $\tau_+$-closed if and only if 

$$\forall a: \N\rightarrow C:  C \supset L_+ (a) := \{ P \in IP(X) \big\vert I^-(\liminf (a) ) = I^-(\limsup (a)) = I^-(P) \} .$$

 $\tau_+$ is manifestly metrizable, finer than $\tau_-$, and shares with $\tau_-$ the property of $i_X$ being a homeomorphism onto its dense image. If $X$ admits a conformal future-compact extension, it is homeomorphic to $IP(X)$ with either topology. The relation between the convergence structures is interesting: For each $p \in X$, each sequence $a: \N\rightarrow J^+ (p) $ has a subsequence $b= a \ci j$ convergent in $(C(X), d_1)$ to some $A \in C(X)$ with $I^-(A) \subset A$ which is, in general, decomposable. Let $U$ be a maximal indecomposable past subset in $A$. Then $b (n) \rightarrow^{\tau_-}_{n \rightarrow \infty} U$. For details on this, see \cite{oM:TFC}, which also describes $\tau_+$ on $IP(X)$ for multiply warped chronological\ spaces $X$.
In a remarkable application, Costa e Silva, Flores and Herrera \cite{CFH} redefined the notion of black hole using $\tau_+$ and reproduced under technical assumptions (but without ever assuming the existence of a conformal future-compact extension) the classical fact that an outer trapped surface is contained in the complement of the past of those sets $A \in IP(X)$ that are the pasts of future null geodesics of infinite affine length. 

\newpage

\section{From ordered measure spaces to Lorentzian length spaces}

Lorentzian length spaces are a facinating approach to "synthesize" Lorentzian geometry. In their paper \cite{mKcS} (Prop. 5.8), when combined with \cite{pCjG}, Cor.1.17, Kunzinger and S\"amann show that the assignment $ \underline{\KS} : (M,g) \mapsto (M, \ll, \leq , \sigma)$ (where $\sigma$ is the signed Lorentzian distance w.r.t. $g$) {\em almost} takes the category of Lipschitz continuous spacetimes into the category of Lorentzian pre-length spaces. The only missing piece is that there is no natural metric --- actually, there {\em cannot} be a natural metric (see Eq. \ref{NOGO}). A {\bf Lorentzian pre-length space} (cf \cite{mKcS} and the next subsection) is a $5$-tuple $(X, d, \leq, \ll, \sigma)$ where $(X,d)$ is a metric space, $(X,\leq)$ is a partially quasi-ordered space, $\ll \subset \leq$ is transitive, and $\sigma: X \times X \rightarrow [0; \infty]$ is lower semi-continuous and satisfies

\begin{enumerate}
\item $\forall x,y,z \in X: x \leq y \leq z \Rightarrow \s (x,z) \geq \s (x,y ) + \s (y,z)$ (conditional inverse triangle inequality),
\item $\forall x,y \in X: \s (x,y ) >0 \Leftrightarrow x \ll y $ (consistency).
\end{enumerate} 

A {\bf Lorentzian length space} is then defined by the three additional requirements causal path-connectedness, local causal closedness, and localizability as in Def. 3.22 of \cite{mKcS}, and finally the property $K^- \ci L^- (\s) = \s $ where the map $L^- $ maps a Lorentzian distance $\s$ to a Lorentzian length functional $L^-(\s)$ for (rectifiable, e.g. locally Lipschitz) causal curves $c: I \rightarrow X$ via  

$$ L^-(\s) (c) =  \inf \{ \sum \s (c(P_k), c(P_{k+1})) \vert P {\rm \ partition \  of \ } I  \}$$

and $K^-$ assigns to a Lorentzian length functional $l$ a Lorentzian distance $K^-(l)$ via 

$$ K^-(l) (p,q) := \sup \{ \{ l(c) \vert c: p \leadsto q {\rm \ causal \ curve}\} \cup \{ 0 \} \} ,$$

and as before we write $c: x \leadsto y$ for a curve $c$ from $x$ to $y$. The Riemannian counterparts $K^+, L^+$ of $K^-, L^-$ satisfy $ K^+ \ci L^+ \geq Id$ and $ L^+ \ci K^+  = Id$. We also have $ K^- \ci L^- \leq Id$ and $ L^- \ci K^-  = Id$, which can be shown by mimicking the proof in \cite{BBI}, 2.3.12 {\em mutatis mutandis}, i.e. inverting the order, replacing suprema with infima etc. A Lorentzian length space is called {\bf globally hyperbolic (g.h.)} iff it is {\bf non-totally imprisoning} (i.e., for each $K \subset X$ compact there is a finite upper bound for the $d$-length of curves in $K$) and all causal diamonds $J(p,q) := J^+(p) \cap J^-(q)$ are compact. Let $\CGH$ be the category of g.h. spacetimes (cf. Sec.1). Subsections 3.1-3.3 contribute the following: 

\begin{enumerate}
\item We show that, despite of the nonexistence of a functor from $\CGH$ to the category of Lorentzian pre-length spaces by the lack of a metric, there {\em is} such a functor for a slightly weakened target category where we only require metrics on a large family of subsets. 
\item Secondly, we suggest a data reduction defining Lorentzian length spaces from the sole datum of a signed length function $\s$ and explain how to recover the other data if $\s$ is $\tau_+$-continuous.
\item Finally, we show that there is a functor from a category of ordered measure spaces to the category of Lorentzian length spaces that is inverse to the functor found recently in $\cite{rMcS} $ by McCann and S\"amann if restricted to the image of the functor of Item 1.
\end{enumerate}

\subsection{Weak Lorentzian length spaces}

In Sec. \ref{LorFin}, we saw that whereas there is a functor from the category of Cauchy slabs to a category of metric spaces, there is no such functor from the entire category of globally hyperbolic manifolds. Still we would like to embed the latter in a natural way into a synthesized approach in the spirit of Lorentzian length spaces. Thus, let us weaken a bit the definition and call a tuple $(X,F, \leq, \ll , \s) $ {\bf weak Lorentzian pre-length space} (whose category we denote by $\CWP$) iff $X, \leq, \ll, \s$ are as above and $F$ is an object as follows: Let $OC(X)$ be the set of the set of closed subsets $A$ of $X$ that are future and past compact (i.e. $J^\pm (p,A) $ compact for all p$ \in A$), partially ordered by inclusion. We call a subset $A$ of $X$ {\bf uniformly full} iff there is $\e>0$ such that for all $a \in A$ there is $b \in A$ with $\vert \s (a,b) \vert \geq \e$. In other words, if we define the thickness $T(X)$ of $X$ as $T(X) := \inf \{ \sup \{ \vert \s (a,b) \vert : b \in X\} : a \in X \} $, then $X$ is uniformly full iff $T(X) >0$, which is satisfied for $\{ (x,y) \in \R^{1,1} : \vert x \vert <1\} $ but not for $\{ (x,y) \in \R^{1,1} : x < \max \{ 1, 1/y\} \} $. Let $OC^u (X)$ the subset of closed uniformly full elements of $OC(X)$. Let $PM(X)$ be the set of pseudo-metrics on $X$ (i.e., symmetric maps $X \times X \rightarrow [0; \infty) $ satisfying the triangle inequality), partially ordered by pointwise comparison. A {\bf local metric on $X$} is a monotonously increasing map from $OC(X)$ to $PM(X)$ s.t. for each $K \in OC^u(X)$, $f(K) \vert_{K \times K}$ is a complete metric on $K $. We require $F \vert_{OC(X)}$ to be a local metric. We will see that a local metric still suffices for the proof of the limit curve theorem\footnote{it would even suffice to specify only a local Lipschitz structure, mapping each $A \in OC(X)$ to a Lipschitz equivalence class of metrics on $A$. However, in order to make the definition functorial while keeping it as close as possible to the original definition, we require the choice of metrics instead of Lipschitz classes and also make the family of appropriate subsets as large as possible.}. At the same time, any complete metric $d$ on $X$ gives rise to a local metric $F: A \mapsto d \ \forall A \subset X$. A {\bf causal curve} in a pre-length space $X$ is a continuous monotonously increasing map $c$ from a real interval $(I, \leq)$ to $(X, \leq)$ that is {\bf locally Lipschitz} in the following sense: There is $ n \in \N \cup \{\infty\}$, a locally finite open covering $\{ U_i \vert i \in \N_n \} $ of $c(I)$ and $ V_i \in OC(X) $ for each $i \in \N_n$ such that $\cl (U_i) \subset V_i$ and $c \vert_{c^{-1} (U_i)}$ is Lipschitz w.r.t. the local metric $d_{U_i}$.

A {\bf weak Lorentzian length space} is a weak Lorentzian pre-length space with $K^- (L^- (\s)) = \s$ (implying causal path-connectedness) that is moreover {\bf locally causally closed} (i.e. every point has a neighborhood $U$ with $\leq \cap (U \times U )$ closed in $U \times U$), and {\bf localizable} (i.e. each point $x \in X$ has a causally convex neighborhood $U_x$ admitting maximizing curves in $U_x$ for each two causally related $p,q \in U_x$ and with $ \sup \{ \ell^d (c) \vert c: I \rightarrow U_x {\rm \ causal \ curve}\} < \infty$) where $d = d_{V}$ for some $V \in OC(X)$ containing $\cl (U_x)$. Localizability implies $\ell (c) \in (0; \infty)$ for each timelike curve $c$ defined on a compact interval, and together with the fact that $\tau_+$ implies local causal convexity, localizability also implies that the space is non-totally imprisoning. We denote the category of weak Lorentzian length spaces by $\CW$.

Let us point out the differences to the definition of Lorentzian length spaces as in \cite{mKcS}: Apart from replacing the datum of local metrics instead of a metric already in the definition of a weak pre-length space, we also renounce the second requirement in the original definition of localizability, which is a bit unappropriate for our aims as it would forbid future or past boundaries --- instead, we can restrict to the subset of "full" points $p$, i.e., those with $I^\pm (p) \neq \emptyset$, in certain applications.

The only conditions that involve $d$ are first causal path-connectedness with the requirement of Lipschitzness of causal curves, and second localizability --- note that both are well-defined for strong equivalence classes. We denote the category of strongly causal $C^2$ Lorentzian manifolds-with-boundary and causal embeddings with causally convex images by $\CCS$. We will see that there is a functor from $\CCS $ to $\CW$, as we can construct {\em local} metrics in a functorial way. But let us postpone this result a bit, as it is most transparent after a further data reduction step in the next subsection. 

%
%

The main instance where the metric in the definition of Lorentzian length spaces is used in $\cite{mKcS}$ is the proof of the Limit Curve Theorem, in order to define uniform Lipschitzness for causal curves in a fixed causal diamond. As the most important applications of this theorem are in compact sets, we can safely replace a metric $d$ by any other metric $d'$ strongly equivalent to $d$ on each compact set. Indeed, if we have a sequence of causal curves $c_i: p \leadsto q$, and if $X$ is a globally hyperbolic $C^1$ spacetime, $ J(p,q)$ is compact and carries a canonical Lipschitz structure, as $C^1$ maps are locally Lipschitz. But if $X$ is more general (a $C^0$ spacetime, a general Lorentzian length space), this is not true any more in general, as there are compact metric spaces $(X,d)$ whose topology can also be induced by a metric $d'$ not strongly equivalent to $d$. In this new setting, we can repeat the proof of the limit curve theorem Th. 3.7 in \cite{mKcS} (for the case that all involved curves take values in a fixed compact subset) for any globally hyperbolic {\em weak} Lorentzian length space and also conclude:

\begin{Theorem}
Let $X$ be a globally hyperbolic weak Lorentzian length space. Then $X$ is strongly causal, the Lorentzian distance function $\s$ on $X$ is finite and continuous, and $X$ is geodesic, i.e. the supremum in the definition of $L^-$ is attained. 
\end{Theorem}

\V{\bf Proof} verbatim as in the proofs of Th.3.26 and Th.3.28 in \cite{mKcS}, replacing $d$ with $d_{J(p,q)}$ \hfill \qed

\subsection{Almost Lorentzian length spaces}

Even more, let us try to define an almost Lorentzian length space by the lone datum of the distance function, a step of radical data reduction, whose justification we postpone a bit:

\begin{Definition}
An {\bf almost Lorentzian pre-length space} is a tuple $(X, \s)$ where $X$ is a set and $\s: X \times X \rightarrow \R$ is antisymmetric and satisfies the conditional inverse triangle inequality $\s(x,z) \geq \s(x,y) + \s (y,z)$ whenever $\s(x,y) , \s (y,z) > 0$, and is continuous w.r.t. $\tau_+(\a (\ll)) $. \\ We denote the category of almost Lorentzian pre-length spaces by $\CALP$ (the morphisms being those bijections $F:X \rightarrow Y$ pulling back the functions, i.e. $\s_Y \circ (F \times F) = \s_X$). An object $X$ of $\CALP$ is called {\bf globally hyperbolic} iff for each $p,q \in X$, the causal diamond $J(p,q)$ is compact.
\end{Definition}

The assignment $ {\rm ldf}: \CCS \rightarrow \CALP $, $(M,g) \mapsto (M, \s_g)$ is obviously a functor (but it should be kept in mind that not every Lorentzian isometric time-oriented embedding is resp. induces a morphism in $\CPOM$ resp. $\CALP$ but only those with causally convex images).

\bigskip

\V Now we want to define a functor ${\rm drf}: {\rm }\CALP \rightarrow \CWP$. First, for $(X, \s) \in {\rm Obj} (\CALP)$, we {\em define} a chronological relation $\ll$ by the second item above (consistency), and $\leq := \a (\ll) \supset \ll $ as in the previous section, and induce the topology $\tau_+$. For each $K \subset X$ define a generalized pseudometric $d_K$ on $X$ by $d_K(p,q) := \sup \{ \vert \s(p, c) - \s (q, c) \vert : c \in K \} $, thus $d_X(p,q) = d_{{\rm Id}, \infty}$. Finally, we define ${\rm fks} := {\rm drf} \ci {\rm ldf}$, which obviously extends the functor $\underline{{\rm fks}}$ from the beginning of Sec.3 by the datum of a local metric. 

\begin{Theorem}
${\rm drf}: (X, \s) \mapsto (X, K \mapsto d_K, \leq, \ll, \s)  $ as above is a functor $\CALP \rightarrow \CWP$. \\ Every globally hyperbolic object of $\CALP$ is mapped to a globally hyperbolic object in $\CWP$.\footnote{The second assertion is nontrivial: Global hyperbolicity in $\CALP$ was defined without non-imprisonment.}
\end{Theorem}

{\bf Proof.} The Noldus metric $d_K \vert_{K \times K}$ is complete for each $K \in OC^u(X)$ as by uniform fullness of $K$, for each Cauchy sequence $a$ in $K$, the cones $J^\pm (a(n))$ eventually intersect each other (otherwise $d_K (a(n), a(m)) > \e$), thus, by future and past precompactness of $K$, $a$ is eventually contained in some compact subset of $X$, and, consequently, contains a convergent subsequence, which implies that $a$, being Cauchy, converges itself, and, by closedness of $K$, to a point in $K$. 

Let $(X, \s)$ be a globally hyperbolic object of $\CALP$, let $A \subset X $ be compact, then there are finite sets $P:= \{ p_i \vert i \in \N_n \}$ and $Q:=  \{ q_i \vert i \in \N_m \}$ such that $A \in I^+ (P) \cap I^-(Q)$. Let us define a time function $t$ on $A$ by $\forall x \in A: t(x) := \sum_{j=0}^n \s (p_j, x) + \sum_{k=0}^m \s (x, q_k)$, then the conditional inverse triangle inequality implies for $y \geq x \geq p \geq p_i $ that $\vert \s (p,y ) - \s (p,x) \vert \leq \vert \s (p_i , y) - \s (p_i, x) \vert \leq t(y) - t(x) $, so $D := \sup \{ t(y) - t(x) \vert x,y \in A\} $ is an upper bound for the length of each causal curve in $A$. \hfill \qed

\bigskip

Finally, we define an {\bf almost Lorentzian length space} to be an object of ${\rm drf}^{-1} (\CW)$, and we denote the category of almost Lorentzian length spaces by $\CALL$.

\begin{Theorem}
$\KS$ takes values in $\CW$, thus ${\rm ldf}$ takes values in $\CALL$.
\end{Theorem}

{\bf Proof.} The topology $\tau_+$ applied in fks recovers the manifold (Alexandrov)  topology, which is locally compact. In \cite{mKcS}, Ex.3.24 (i), it has been shown that the two requirements of path connectedness via Lipschitz curves and of localizability are satisfied for an arbitrary Riemannian metric on the spacetime. As both requirements concern compact sets and are well-defined for strong equivalence classes, it is enough to show that on compact subsets the Noldus metric is strongly (i.e., Lipschitz) equivalent to any Riemannian metric, which follows from the first displayed equation of Theorem 5 in \cite{oM:LGH}. \hfill \qed

\bigskip

\V {\bf Remark.} It is instructive to see localizability in the framework of Sec. \ref{LorFin} for the Riemannian metric induced by $\Phi_{g,f,2}$:
For $ \s (q,p_n^+) , \s (p_n^-, q)   \in (0; 1/n)$, the extrinsic $d_{J(p_n^-, p_n^+)}$-diameter $d$ of $J(p_n^-, p_n^+)$ (coming from an embedding $\Phi_g$ into $L^2(X)$ as in the first section) tends to zero for $n \rightarrow \infty$. The norm of the second fundamental form on every $\Phi_g (A)$ is monotonically increasing w.r.t. inclusion, so for a small enough neighborhood $U$ the condition $\vert \vert S_M^H \vert \vert < 4 \sqrt{2} \cdot (3 \delta )^{-1}$ is satisfied, with $\delta := \diam_{L^2 (X)} (U)$. Then, by Theorem \ref{Hilbert}, for the intrinsic diameter $D$ we obtain $D < 3 \delta $, that is, the Riemannian length of causal curves in $U$ is bounded above by $ 3 \diam_{L^2 (X)} (U) $.

\subsection{Natural Lorentzian distances from orders and measures}

Keeping in mind Berger's quote above, we see that there is actually a Lorentzian analogue to metric measure spaces: (partially) ordered measure spaces, forming the category $\CPOM$. Firstly, the map $\SEP: (M,g) \mapsto (M, \leq, \Vol_g)$ is an obvious (injective, see \cite{oM:LGH}, proof of Th. 1 (ii)) functor from the category $\CCS$ to $\CPOM$, which has many more objects than $\SEP (\CCS)$ --- e.g., a short consideration shows that the multiply warped chronological spaces from \cite{oM:TFC} carry not only a natural causal structure  but also a natural Borel measure if the slice is of constant finite Hausdorff dimension. Ordered measure spaces improve one aspect of metric measure spaces in that they assign two degrees of freedom to two independent properties of physical interest each: on one hand the causal structure, encoding the mutual  dependence between physical phenomena, and on the other hand the volume of the spacetime regions.
Very likely, the categories $\CPOM$ and $\CALP$ will turn out to be of considerable physical interest as they harbour many possible degenerations of causal spacetimes and still could permit the formulation of physical theories.

McCann and S\"amann \cite{rMcS}, along the central and very fruitful idea to replace, in the definition of Hausdorff dimension and measure, the balls with causal diamonds, construct a natural measure on Lorentzian pre-length spaces (and thereby a functor $\MCS: \CWP \rightarrow \CPOM $), equivalent to the following: For $A \subset X$, $N>0$ we define (keeping in mind $J(x,y) \in OC(X) \ \forall x,y \in X$):

$CC_{\delta} (A) := \{ (p,q) \in (X^{\N})^2 \vert A \subset \bigcup_{k=1}^{\infty} J(p(k), q(k)) \land \diam (J (p(k), q(k))) < \delta \forall k \in \N  \},$

(where $\diam$ is w.r.t. the metric $d_{U_k}$ for $U_k:= J(p_k,q_k)$),

$$ \la_N ((p,q)) := \omega (N) \cdot \sum_{k=1}^\infty \s (p(k),q(k))^N \forall (p,q) \in CC_{\delta} (A), \ {\rm where} \  \omega(N):= \frac{\pi^{\frac{N-1}{2}}}{N \cdot \Gamma (\frac{N+1}{2}) \cdot 2^{N-1}}, $$

$$ \mu_{N, \delta} (A) := \inf \{ \la_N(D) \vert D \in CC_{\delta} (A)\}, \qquad  \mu_N (A) := \lim_{\delta \rightarrow 0} \mu_{N, \delta} (A),  $$

\V The outer measure $\mu_N$ induces a unique measure on the Borel subsets. The original definition in \cite{rMcS} deviates from the above insofar as the authors require $\diam_d (J(p(k),q(k))) < \delta$, i.e., the diameter is computed by means of the unique metric $d$, which is part of the equipment of the original definition of pre-length space but not of almost pre-length space (and again, the choice of a local strong equivalence class instead of a local metric would suffice here).


\bigskip

We want to take another step in transferring physical theories from $\CCS$ to $\CPOM$ or to $\CALP$, showing that there is a functor $\F: \CPOM \rightarrow \CALP $ with $\KS= \F \ci \SEP$, reconstructing the Lorentzian distance function from order and measure alone, even when no tangent space is available. 

\bigskip

The volume of a ball of radius $R$ in Euclidean $\R^n$ is ($\Gamma$ being Euler's gamma function)

$$V_n(R) = \frac{\pi^{n/2}}{\Gamma (\frac{n}{2} +1)} R^n ,$$

and the volume $K_n(h)$ of the $n$-dimensional cone of height $h$ and aperture $\pi/2$ over a ball is

$$K_n(h)  = \frac{\pi^{\frac{n-1}{2}}}{ n \cdot \Gamma(\frac{n-1}{2} +1)} h^n, {\rm \ thus, \ solved \ for \ } h : \   h = \sqrt[n]{\frac{n \Gamma (\frac{n-1}{2} +1)}{\pi^{\frac{n-1}{2}}} \cdot K_n(h)} . $$

On the other hand, we see that $K_n(h/2) = 2^{-n} K_n(h)$. Define (for $a,c \in M$ with $\Vol (J(a,c)) \neq 0$)

$$ \Phi(a,c)  :=  \Phi_{\leq, \Vol} (a,c) :=  \sup \Big\{ \frac{\Vol(J(a,b) ) + \Vol(J(b,c) )}{\Vol (J(a,c)) } \Big\vert b \in J(a,c)  \land \Vol (J(a,b)) = \Vol (J(b,c))  \Big\} \in (0; 1]$$

where the last estimate is due to $J(a,b) \cup J(b,c) \subset J(a,c) \ \forall b \in J(a,c)$. In $\R^{1,n}$, the supremum is attained at the midpoint $b$ of the unique geodesic from $a$ to $c$, and $\Vol^{1/n}(J(a,c) )= \Vol^{1/n} (J(a,b)) + \Vol^ {1/n}(J(b,c))$.
With this in mind, for an ordered measure space $X \ni b$ we define the {\bf dimension} $\Dm (b)$ {\bf of $X$ at $b$} (which recovers the dimension $n $ at every $b \in X$ if $X$ is a spacetime) by

$$ R (b)  := \{ (p,q) \in J^-(b) \times J^+(b) \vert \Vol (J(p,q) ) \neq 0 \} , \ \  \Dm(b) := \limsup_{R(b) \ni (p,q) \rightarrow (b,b)} (-\log_2 (\Phi (a,c)) +1) \in [1; \infty].$$

Then we define the Lorentzian length $\ell(c) \in [0; \infty]$ of a timelike curve $c: I \rightarrow X$ by

$$ \ell(c) = \inf \Big\{ \sum_{k=1}^N \sqrt[\Dm(p_k)]{\frac{\Dm(p_k) \Gamma (\frac{\Dm(p_k)-1}{2} +1)}{\pi^{\frac{\Dm(p_k)-1}{2}}} \cdot \Vol (J(p_k, p_{k+1}))}  \Big\vert \{t_0 ,... t_N \} {\rm partition \ of \ } I,  p(n) := c(t_n) \Big\}. $$

Finally ${\rm fld} (X,\leq, \Vol) := (X, K^-(\ell))$. Geodesic normal coordinates (in which the first derivative of the metric vanishes at $0$) show that on $\SEP (\CCS)$, the so defined Lorentzian distance is indeeed $\s_g$:

\begin{Theorem}
For the functors $\F$, $\SEP$, $\KS$ defined as above we have $\F \circ \SEP = \KS $. 

\end{Theorem}

\V{\bf Proof.} As $g$ is a differentiable metric on the spacetime $M$, we can define Lorentzian normal coordinates at a point $p$. Those fix a Euclidean metric $E$ on $T_p M$. As $d_pg = 0 $ in normal coordinates, the one-parameter family of rescaled metrics $g^\e := \frac{1}{\e} \cdot g \vert_{B^E(0, \e)}  $ converges with $\e \rightarrow 0$ to the Minkowski metric in the $C^0$ norm. For $ \exp_p (u) = a \ll p \ll c = \exp_p (v)$, the ratio $\Phi_g (a,c) := \Phi_{\Vol_g, \leq_g} $ is invariant under scaling $g \leadsto \la g$ of $g$ and continuous in $g$, thus $\Phi_{g^{1/m}} (\exp_p (\frac{1}{m}u),\exp_p (\frac{1}{m}v)) \rightarrow_{m \rightarrow \infty} \Phi_{g_{1,n}} (\tilde{a},\tilde{c}) = 2^{-n}$, where $\tilde{a} $ and $\tilde{c}$ are two points of Lorentzian distance $1$ in $\R^{1,n}$, therefore the dimension is recovered correctly (and the $\limsup$ is a true limit here). The ratio $\beta_g (a,c):= \frac{(\Vol_g (J(a,c)))^{1/n}}{\s_g (a,c)} $ is invariant under scaling $g \leadsto \la g$ of $g$ and continuous in $g$ as well, thus in the formula for $\ell$, we get that $(\Vol_g (J(c(t_k), c(t_{k+1}))))^{1/n}/\s_g(c(t_k, c(t_{k+1})))$ converges to $1$ uniformly for all $k$ as the fineness of the partition goes to $0$, which implies $\ell := L(\s_g)$. \hfill \qed

\newpage

McCann and S\"amann, in \cite{rMcS}, show $ \MCS \ci \KS= \SEP $, thus on the respective images of $\CCS$, the maps $\F$ and $\MCS$ are inverse to each other, thus injectiviy of sep implies injectivity of mcs.

{\bf Question:} On which subcategories of $\CPOM$ resp. $\CALP$ are $\MCS$ and $\F$ inverse to each other? 

A short consideration shows that on multiply warped chronological spaces this is true. In general, constancy of the Lorentzian Hausdorff dimension seems a good condition in this context. 

\bigskip

Here is a panoramic diagram containing most of the functors discussed so far. The vertical functors of its ground floor are the obvious forgetful ones, the vertical lines of its upper floor are the involved synthesization functors. e.g. the functor ${\rm dst}$ assigning to each object of $C_n^{+,c} $, which is a Riemannian manifold-with-boundary, the corresponding metric space. Dashed lines indicate that the respective functor is defined only in a (considerably large) subset. The polygons of the diagram are commutative under the inscribed conditions. Furthermore, ${\rm tau} : {\bf POS} \rightarrow \CTop$ is the functor $(X, \leq) \mapsto (X, \tau_+(\beta (\leq)))$. The functor assigning to an almost Lorentzian length space $X$ or to an ordered measure space the metric space $(X, K^- \ci L^- (\Phi_f^*(d_{L^p(X)}))) $ is called ${\rm phi} (p,f)$, and ${\rm phl} (p,f) (X, g) = (X, d\Phi^* (\langle \cdot , \cdot  \rangle_{L^2 (X)})) $. We denote the subcategory of ${\bf C}^-_n$ consisting of the uniformly full objects by ${\bf C}^{-,u}_n$. For $f \in {\rm Adm} \cap C^4 (\R, \R)$ vanishing to fourth order at $0$ we get 

\medskip
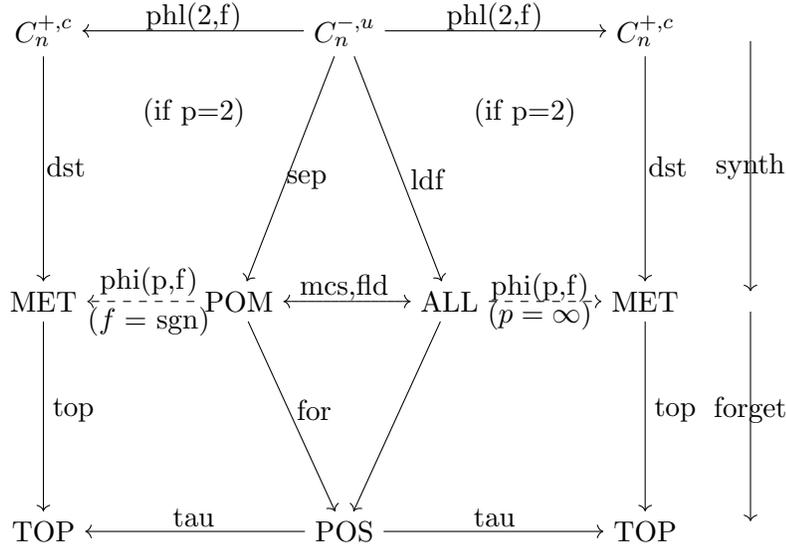
\begin{center}
\tikzset{every node/.style={align=center}} 

\begin{tikzpicture}[x=2cm,y=1.8cm, scale=1.0, every node/.style={transform shape}]
\node at (0,0) (n) {$C_n^{+,c}$};
\node at (0,-2) (m) {MET};
\node at (0.2,-2.8) (tr) {top};
\node at (0,-3.7) (p) {TOP};
\node at (1.3,-2) (c) {POM};
\node at (0.7,-1.85) (cf) {phi(p,f)};
\node at (0.7,-2.15) (cu) {($f= $ sgn)};
\node at (2,-3.7) (e) {POS};
\node at (2.7,-2) (f) {ALL};
\node at (2,0) (g) {$C_n^{-,u}$};
\node at (4,-3.7) (i) {TOP};
\node at (4,-2) (j) {MET};
\node at (4,0) (y) {$C_n^{+,c} $};
\node at (4.7,0) (u) {};
\node at (4.7,-1) (k) {synth};
\node at (4.7,-2) (q) {};
\node at (4.7,-2.8) (r) {forget};
\node at (4.7,-3.7) (h) {};
\node at (2,-1.9) (a) {mcs,fld};
\node at (1.75, -1.1) (s) {sep};
\node at (2.55, -1.1) (t) {ldf};
\node at (3.3,-1.9) (a) {phi(p,f)};
\node at (3.3,-2.1) (as) {($p= \infty$)};
\node at (3,-3.6) (v) {tau};
\node at (1.8,-2.8) (w) {for};
\node at (4.2,-2.8) (x) {top};
\node at (3,0.1) (z) {phl(2,f)};
\node at (1,0.1) (z) {phl(2,f)};
\node at (4.15,-1) (z) {dst};
\node at (0.15,-1) (z) {dst};
\node at (3.2,-0.6) (z) {(if p=2)};
\node at (1,-0.6) (z) {(if p=2)};
\node at (1,-3.6) (zu) {tau};
\draw [->] (g)--(n);
\draw [->] (u)--(q);
\draw [->] (n)--(m);
\draw [dashed, ->] (c)--(m);
\draw [->] (c)--(e);
\draw [->] (g)--(y);
\draw [->] (y)--(j);
\draw [->] (g)--(f);
\draw [->] (g)--(c);
\draw [->] (j)--(i);
\draw [dashed, ->] (f)--(j);
\draw [dashed, ->] (c)--(f);
\draw [dashed, ->] (f)--(c);
\draw [->] (e)--(i);
\draw [->] (e)--(p);
\draw [->] (m)--(p);
\draw [->] (f)--(e);
\draw [->] (q)--(h);
\end{tikzpicture}

\end{center}

\medskip

It is conceivable that ordered measure spaces and almost Lorentzian pre-length spaces form morally equivalent categories, the latter comprising everything in the single datum of Lorentzian length, the former separating the aspects of influenceability and measure/overall importance of events.

\bigskip

{\bf Data availability statement:} No experimental data have been produced for this article.

\end{document}